\documentclass[11pt]{amsart}
\usepackage{latexsym,amssymb,amsfonts,amsmath,graphicx,fullpage,url}
\usepackage[vcentermath,enableskew]{youngtab}
\usepackage{color}
\usepackage{cite}
\usepackage[all]{xy}
\usepackage{verbatim}
\usepackage{tikz,hyperref}
\usepackage{caption}
\usepackage{subcaption}
\newtheorem{theorem}{Theorem}[section]

\newtheorem{proposition}{Proposition}[section]
\newtheorem{lemma}[theorem]{Lemma}
\newtheorem{corollary}[theorem]{Corollary}
\theoremstyle{definition}
\newtheorem{definition}[theorem]{Definition}
\newtheorem{example}[theorem]{Example}

\theoremstyle{remark}

\numberwithin{equation}{section}

\def\C{{\mathcal{C}}}
\def\S{{\mathcal{S}}}

\def\t{{\mathcal{S}(\beta)}}
\def\P{{\mathcal{P}}}

\newcommand{\be}{\begin{equation}}
\newcommand{\ee}{\end{equation}}

\newcommand{\bd}{\begin{definition}}
\newcommand{\ed}{\end{definition}}
\newcommand{\bt}{\begin{theorem}}
\newcommand{\et}{\end{theorem}}
\newcommand{\bl}{\begin{lemma}}
\newcommand{\el}{\end{lemma}}
\newcommand{\bp}{\begin{proposition}}
\newcommand{\ep}{\end{proposition}}
\newcommand{\bc}{\begin{corollary}}
\newcommand{\ec}{\end{corollary}}

\newcommand{\lmp}[1]{\mathrel{\mathop{\longmapsto}^{\mathrm{#1}}}} 
\newcommand{\sizeQ}{m} 

\newtheorem*{theorem1*}{Theorem \ref{thm:main}}

\newcommand{\old}[1]{}

\font\co=lcircle10
\def\petit#1{{\scriptstyle #1}}

\def\jr{\smash{\raise2pt\hbox{\co \rlap{\rlap{\char'005} \char'007}}
               \raise6pt\hbox{\rlap{\vrule height6.5pt}}
               \raise2pt\hbox{\rlap{\hskip4pt \vrule height0.4pt depth0pt
                width7.7pt}}}}
\def\je{\smash{\raise2pt\hbox{\co \rlap{\rlap{\char'005}
                \phantom{\char'007}}}\raise6pt\hbox{\rlap{\vrule height6pt}}}}
\def\+{\smash{\lower2pt\hbox{\rlap{\vrule height14pt}}
                \raise2pt\hbox{\rlap{\hskip-3pt \vrule height.4pt depth0pt
                width14.7pt}}}}
\def\perm#1#2{\hbox{\rlap{$\petit {#1}_{\scriptscriptstyle #2}$}}%
                \phantom{\petit 1}}

\def\textcross{\ \smash{\lower4pt\hbox{\rlap{\hskip4.15pt\vrule height14pt}}
                \raise2.8pt\hbox{\rlap{\hskip-3pt \vrule height.4pt depth0pt
                width14.7pt}}}\hskip12.7pt}
\def\textelbow{\ \hskip.1pt\smash{\raise2.8pt%
                \hbox{\co \hskip 4.15pt\rlap{\rlap{\char'005} \char'007}
                \lower6.8pt\rlap{\vrule height3.5pt}
                \raise3.6pt\rlap{\vrule height3.5pt}}
                \raise2.8pt\hbox{%
                  \rlap{\hskip-7.15pt \vrule height.4pt depth0pt width3.5pt}%
                  \rlap{\hskip4.05pt \vrule height.4pt depth0pt width3.5pt}}}
                \hskip8.7pt}

\title{Subword complexes via triangulations of root polytopes}
\author{Laura Escobar}
\address{Laura Escobar,
 Department of Mathematics, University of Illinois at Urbana-Champaign, Urbana IL 61801  \newline lescobar@illinois.edu
}

\author{Karola M\'esz\'aros}
\address{Karola M\'esz\'aros, Department of Mathematics, Cornell University, Ithaca NY 14853  \newline karola@math.cornell.edu
}
\thanks{M\'esz\'aros was partially supported by a National Science Foundation Grant  (DMS 1501059).}
\date{\today}
\begin{document}

 \maketitle
 
 \begin{abstract}
 Subword complexes are simplicial complexes introduced by Knutson and Miller to illustrate the combinatorics of Schubert polynomials and determinantal ideals. They proved that any subword complex is homeomorphic to a ball or a sphere and asked about their geometric realizations. We show that a family of subword complexes can be realized geometrically via regular triangulations of root polytopes. This implies that a family of $\beta$-Grothendieck polynomials are special cases of reduced forms in the subdivision algebra of root polytopes. We can also write the volume and Ehrhart series of   root polytopes in terms of $\beta$-Grothendieck polynomials. 
 \end{abstract}
 \tableofcontents 

\section{Introduction}
 
 In this paper we provide  geometric realizations of pipe dream complexes $PD(\pi)$ of permutations $\pi=1\pi'$, where $\pi'$ is a dominant permutation on $2,3,\ldots, n$ as well as the subword complexes that are the cores of the pipe dream complexes $PD(\pi)$.  We realize  $PD(\pi)$  as (repeated cones of) regular  triangulations  of the root polytopes $\P(T(\pi))$.

 Since the appearance of Knutson's and Miller's work in \cite{subword, annals} there has been a flurry of research into the geometric realization of subword complexes with progress in realizing  families of spherical subword complexes \cite{stump,cesar,assoc,abrick,SerranoStump,subwordcluster,reviewedfansub}.  This paper is the first to succeed in realizing a family of subword complexes which are homeomorphic to balls.  
 
Subword complexes were first shown to relate to triangulations of root polytopes by M\'esz\'aros in \cite{groth}, where the author gives a geometric realization of the pipe dream complex of $[1,n,n-1,\ldots,1]$ and whose work served as the stepping stone for the present  project. In the papers \cite{root1, prod, h-poly1, h-poly2}  M\'esz\'aros  studied  triangulations of root polytopes that  we utilize in this work (some of the mentioned papers are in the language of flow polytopes, but in view of \cite[Section 4]{groth} some of their content can also be understood in the language of root polytopes). 
 
 The main theorem of this paper is the following, which has several interesting consequences explored in the paper. For the definitions needed for this theorem see the later sections.
 
 \bt  \label{thm:intro} Let $\pi=1\pi' \in S_n$, where $\pi'$ is dominant. Let $\mathcal{C}^2(\pi)$ be the core of $PD(\pi)$ coned over twice. The   canonical triangulation of the root polytope $\P(T(\pi))$ (which is a regular triangulation) is a geometric realization of $\mathcal{C}^2(\pi)$.
 \et 
 
 The roadmap of this paper is as follows. In Sections \ref{sec:pipe} and \ref{sec:root} we explain the necessary background about subword complexes and root polytopes, respectively. In Section \ref{sec:constr} we prove a geometric realizations of pipe dream complexes $PD(\pi)$ of permutations $\pi=1\pi'$, where $\pi'$ is a dominant permutation on $2,3,\ldots, n$, via triangulations of   root polytopes $\P(T(\pi))$. In Section \ref{sec:groth} we use the previous result to show that $\beta$-Grothendieck polynomials are special cases of reduced forms in the subdivision algebra of root polytopes while in Section \ref{sec:cor} we show how to express the volume and Ehrhart series of root polytopes in terms of Grothendieck polynomials. Section \ref{sec:unique} is devoted to proving a certain uniqueness property of the reduced form in the subdivision algebra that we used in Section \ref{sec:groth}.
 
\section{Background on pipe dream complexes}
\label{sec:pipe}

We let $S_n$ denote the set of permutations of size $n$.
\bd The {\bf (Rothe) diagram} of a permutation $\pi\in S_n$ is the collection of boxes $D(\pi)=\{(\pi_j,i): i<j, \pi_i>\pi_j\}$. It can be visualized by considering the boxes left in the $n\times n$ grid after we cross out the boxes appearing south and east  of each 1 in the permutation matrix for $\pi$.
\ed

\begin{figure}[h]
\begin{tikzpicture}
\draw (0,0)--(2,0)--(2,2)--(0,2)--(0,0);
\draw[]
    (2,1.25)
 -- (1.75,1.25) node {$\bullet$}
 -- (1.75,0);
\draw[]
    (2,1.75)
 -- (0.75,1.75) node {$\bullet$}
 -- (0.75,0);
 \draw[]
    (2,.75)
 -- (1.25,.75) node {$\bullet$}
 -- (1.25,0); 
 \draw[]
    (2,.25)
 -- (.25,.25) node {$\bullet$}
 -- (.25,0);
 
 \draw (0,.5)--(.5,.5)--(.5,2);
 \draw (0,1)--(.5,1);
 \draw (0,1.5)--(.5,1.5);
 
\draw (1,1)--(1.5,1)--(1.5,1.5)--(1,1.5)--(1,1);
 
\end{tikzpicture}\caption{The diagram for $\pi=[4132]$.}
\end{figure}
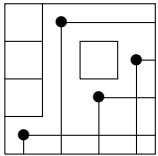


Notice that no two permutations can give the same diagram. We will consider permutations of the form $\pi=1\pi'$ where $\pi'$ is a {\bf dominant permutation} of $\{2,\ldots,n\}$, i.e., the diagram of $\pi$ is a partition with north-west most box at position $(2,2)$. Dominant permutations can be equivalently defined as the $132$-avoiding permutations, and there are Catalan many for fixed size. Our convention is to encode the partition by the number of boxes in each column.

\begin{figure}[h]
\begin{tikzpicture}
\draw (0,0)--(3,0)--(3,3)--(0,3)--(0,0);
\draw[]
    (3,2.75)
 -- (0.25,2.75) node {$\bullet$}
 -- (0.25,0);
\draw[]
    (3,0.25)
 -- (0.75,0.25) node {$\bullet$}
 -- (0.75,0);
 \draw[]
    (3,1.25)
 -- (1.25,1.25) node {$\bullet$}
 -- (1.25,0); 
 \draw[]
    (3,2.25)
 -- (1.75,2.25) node {$\bullet$}
 -- (1.75,0);
 \draw[]
    (3,1.75)
 -- (2,1.75) node {$\bullet$}
 -- (2,0);
 \draw[]
    (3,0.75)
 -- (2.25,0.75) node {$\bullet$}
 -- (2.25,0);
 
 \draw (0.5,0.5)--(1,0.5)--(1,2.5)--(0.5,2.5)--(0.5,0.5);
 \draw (0.5,1)--(1,1);
 \draw (0.5,1.5)--(1.5,1.5);
  \draw (0.5,2)--(1.5,2);
 \draw (1,2.5)--(1.5,2.5)--(1.5,1.5);

\end{tikzpicture}\caption{The diagram for $\pi=[164235]$ which corresponds to $\lambda=(4,2)$}
\end{figure}
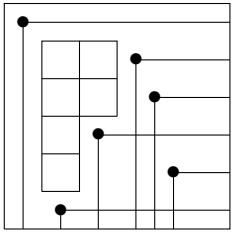

\bd A {\bf pipe dream} for $\pi\in S_n$ is a tiling of an $n\times n$ matrix with two tiles, crosses $\textcross$ and elbows $\textelbow$, such that
\begin{enumerate}
\item all tiles in the weak south-east triangle of the $n\times n$ matrix are elbows, and
\item if we write $1,2,\ldots, n$ on the left and follow the strands (ignoring second crossings among the same strands) they come out on the top and read $\pi$.
\end{enumerate}
A pipe dream is {\bf reduced} if no two strands cross twice.
\ed
\begin{figure}[h]
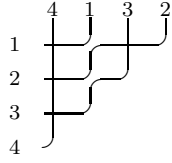
$
\begin{array}{ccccc}
       &\perm4{}&\perm1{}&\perm3{}&\perm2{}\\
\petit1&   \+   &   \jr  &   \+   &  \je   \\
\petit2&   \+ &   \jr  &   \je  &\\
\petit3&   \+  &   \je  &        &\\
\petit4&   \je  &        &        &\\
\end{array}
$
\caption{A reduced pipe dream for $\pi=[4132]$.}\end{figure}

\bd The {\bf pipe dream complex} $PD(\pi)$ of a permutation $\pi\in S_n$ is the simplicial complex with vertices given by entries on the northwest triangle of an $n\times n$-matrix and facets given by the elbow positions in the reduced pipe dreams for $\pi$.
\ed

Pipe dream complexes are a special case of the subword complexes defined by Knutson and Miller in \cite{subword,annals}. We proceed to explain the correspondence. The group $S_n$ is generated by the adjacent transpositions $s_1,\ldots,s_{n-1}$, where $s_i$ transposes $i \leftrightarrow i+1$. Let $Q=(q_1,\ldots,q_\sizeQ)$ be a word in $\{s_1,\ldots,s_{n-1}\}$, i.e., $Q$ is an ordered sequence. A {\bf subword} $J=(r_1,\ldots,r_\sizeQ)$ of $Q$ is a word obtained from $Q$ by replacing some of its letters by $-$. There are a total of $2^{|Q|}$ subwords of $Q$. Given a subword $J$, we denote by $Q\setminus J$ the subword with $k$-th entry equal to $-$ if $r_k\neq -$ and equal to $q_k$ otherwise for, $k=1,\ldots,\sizeQ$. For example, $J=(s_1,-,s_3,-,s_2)$ is a subword of $Q=(s_1,s_2,s_3,s_1,s_2)$ and $Q\setminus J=(-,s_2,-,s_1,-)$. Given a subword $J$ we denote by $\prod J$ the product of the letters in $J$, from left to right, with $-$ behaving as the identity.

\bd \cite{subword,annals}
 Let $Q=(q_1,\ldots,q_\sizeQ)$ be a word in $\{s_1,\ldots,s_{n-1}\}$ and $\pi\in S_n$. The \textbf{subword complex} $\Delta(Q,\pi)$ is the simplicial complex on the vertex set $Q$ whose facets are the subwords $F$ of $Q$ such that the product $\prod(Q\setminus F)$ is a reduced expression for $\pi$.
\ed

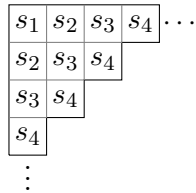
\begin{figure}[h]
\begin{tikzpicture}
\draw (0,0)--(0.5,0)--(0.5,0.5)--(1,0.5)--(1,1)--(1.5,1)--(1.5,1.5)--(2,1.5)--(2,2)--(0,2)--(0,0);
\draw[gray] (0,0.5)--(0.5,0.5)--(0.5,2);
\draw[gray] (0,1)--(1,1)--(1,2);
\draw[gray] (0,1.5)--(1.5,1.5)--(1.5,2);

\node at (0.25,0.25) {$s_4$};
\node at (0.75,0.75) {$s_4$};
\node at (1.25,1.25) {$s_4$};
\node at (1.75,1.75) {$s_4$};
\node at (.25,0.75) {$s_3$};
\node at (0.75,1.25) {$s_3$};
\node at (1.25,1.75) {$s_3$};
\node at (0.25,1.25) {$s_2$};
\node at (0.75,1.75) {$s_2$};
\node at (0.25,1.75) {$s_1$};
\node at (2.3,1.75) {$\cdots$};
\node at (0.25, -0.2) {$\vdots$};

\end{tikzpicture}\caption{Labeling of the entries in the northwest triangle by adjacent transpositions.}\label{fig:pipex}\end{figure}

In this language, $PD(\pi)$ is the subword complex $\Delta(Q,\pi)$ corresponding to the triangular word $Q=(s_{n-1},s_{n-2},s_{n-1},\ldots,s_1,s_2,\ldots,s_{n-1})$ and $\pi$. The correspondence between pipe dreams and subwords is induced by the labeling of the entries in the northwest triangle of an $n\times n$-matrix by adjacent transpositions, as depicted in Figure \ref{fig:pipex}, and by making a $\textcross$ in a pipe dream correspond to a $-$ in a subword and a $\textelbow$ correspond to the $s_i$ in its entry. In order to go from a pipe dream to a subword, we read the entries in the northwest triangle from left to right starting at the bottom.

\bd \label{def:cone} Let cone$(\pi)$ be the set of vertices of $PD(\pi)$ that are in all its facets. We define the {\bf core} of $\pi$, denoted by core$(\pi)$, to be the simplicial complex obtained by restricting $PD(\pi)$ to the set of vertices not in cone$(\pi)$.\ed 

Notice that $PD(\pi)$ is obtained from its core by iteratively coning the simplicial complex core$(\pi)$ over the vertices in cone$(\pi)$. This is a standard definition for simplicial complexes. In the language of pipe dream complexes, the core is the restriction to the entries in the $n\times n$ matrix that are a cross in some reduced pipe dream for $\pi$. Following the correspondence described in Figure \ref{fig:pipex}, this restriction induces a subword $Q'$ of the triangular word and so core$(\pi)$ is the subword complex $\Delta(Q',\pi)$.

\begin{figure}[h]
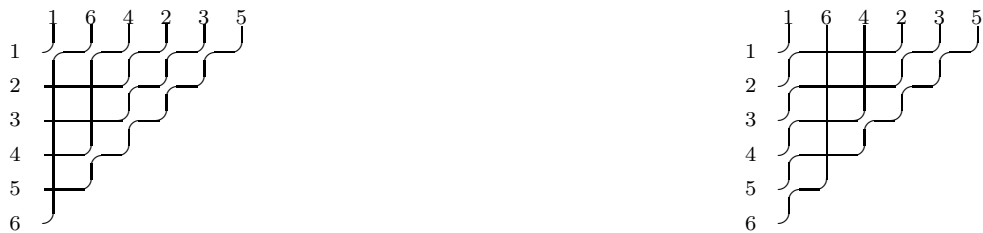

    \centering
    \begin{subfigure}[b]{0.4\textwidth}
        \centering
      $
\begin{array}{ccccccc}
       &\perm1{}&\perm6{}&\perm4{}&\perm2{}&\perm3{}&\perm5{}\\
\petit1&   \jr   &   \jr  &   \jr  &   \jr  &   \jr   &  \je   \\
\petit2&   \+  &   \+ &   \jr &   \jr    &   \je  &\\
\petit3&   \+  &   \+ &   \jr   &    \je    &\\
\petit4&   \+  &    \jr    &  \je      &\\
\petit5&   \+  &     \je   &        &\\
\petit6&   \je  &        &        &\\
\end{array}
$
        \caption{Aligned left is the bottom reduced pipe dream}
    \end{subfigure}
    \hfill
    \begin{subfigure}[b]{0.4\textwidth}
        \centering      $
\begin{array}{ccccccc}
       &\perm1{}&\perm6{}&\perm4{}&\perm2{}&\perm3{}&\perm5{}\\
\petit1&   \jr   &   \+  &   \+  &   \jr  &   \jr   &  \je   \\
\petit2&   \jr  &   \+ &   \+ &   \jr    &   \je  &\\
\petit3&   \jr  &   \+ &   \jr   &    \je    &\\
\petit4&   \jr  &    \+    &  \je      &\\
\petit5&   \jr  &     \je   &        &\\
\petit6&   \je  &        &        &\\
\end{array}
$
        \caption{Aligned up is the top reduced pipe dream}

    \end{subfigure}
    \hfill
    \caption{Two reduced pipe dreams for [164235] obtained by aligning the diagram to the left and to the top.}        \label{fig:topbottomrpd}
\end{figure}

Since we are only considering permutations of the form $1\pi'$ with $\pi'$ a dominant permutation, core$(1\pi')$ is easy to describe. Given a diagram of a permutation there are two natural reduced pipe dreams for $\pi$, referred to as the \textbf{bottom reduced pipe dream of $\pi$} and the  \textbf{top reduced pipe dream of $\pi$}, one obtained by aligning the diagram to the left and replacing the boxes with crosses and the other one by aligning the diagram up. See Figure \ref{fig:topbottomrpd}.

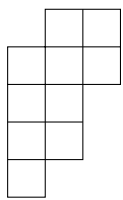
\begin{figure}[h]
\begin{tikzpicture}
 
 \draw (0,1.5)--(1.5,1.5)--(1.5,2.5)--(0.5,2.5)--(0.5,0)--(0,0)--(0,2)--(1.5,2);
 \draw (0,0.5)--(1,0.5)--(1,2.5);
 \draw (0,1)--(1,1);

\end{tikzpicture}\caption{The core region of $[164235]$}\label{fig:crpi}
\end{figure}

The core of $1\pi'$ is the simplicial complex obtained by restricting $PD(\pi)$ to the vertices corresponding to the positions of the crosses in the superimposition of these two pipe dreams. We refer to the region itself as the \textbf{core region}, and denote it by cr$(\pi)$. See Figure \ref{fig:crpi} for an example. Note that different permutations can have the same core region, as is the case for $[15342]$ and $[15432]$.

\begin{figure}[h]
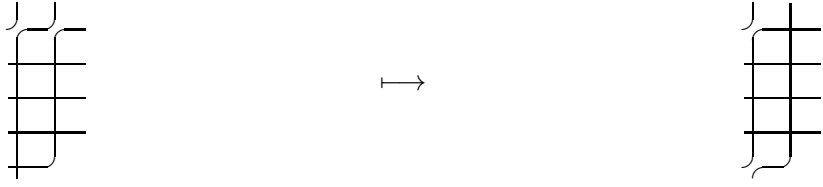

    \centering
    \begin{subfigure}[h]{0.4\textwidth}
        \centering
        
$\begin{array}{cc}
      {\color{white} \perm2{}}&{\color{white}\perm1{}}\\
   \jr   &   \jr   \\
   \+  &   \+ \\
  \+  &   \+\\
   \+  & \+  \\
      \+  & \je  \\
\end{array}$

    \end{subfigure}
    \hfill
    \begin{subfigure}[h]{0.1\textwidth}
        \centering   
        $\displaystyle{\lmp{}}$
        
            \end{subfigure}
    \hfill
    \begin{subfigure}[h]{0.4\textwidth}
        \centering

$\begin{array}{cc}
      {\color{white} \perm2{}}&{\color{white}\perm1{}}\\
   \jr   &   \+   \\
   \+  &   \+ \\
  \+  &   \+\\
   \+  & \+  \\
      \jr  & \je  \\
\end{array}$
        
    \end{subfigure}
    \hfill
    \caption{Ladder move.}\label{fig:ladder}
\end{figure}

In \cite{rc}, Bergeron and Billey introduced an algorithm to construct all reduced pipe dreams for $\pi$. Given a reduced pipe dream $P$ for all permutations $\pi$, a ladder admitting rectangle is a connected $k\times 2$ rectangle inside $P$ such that $k\geq 2$ and the only $\textelbow$ inside this rectangle are in the top row and in the southeast corner, see the diagram on the left in Figure \ref{fig:ladder}. A {\bf ladder move} on $P$ moves the $\textcross$ in the southwest corner of a ladder admitting rectangle to the northeast corner. Notice that the resulting pipe dream is a reduced pipe dream for $\pi$.

\bt\cite{rc}\label{thm:ladder} The set of all reduced pipe dreams of $\pi$ equals the set of pipe dreams that can be derived from the bottom reduced pipe dream by a sequence of ladder moves.
\et

The {\bf boundary} of a pure simplicial complex $\Delta$ is the simplicial complex $\partial\Delta$ with facets the codimension 1 faces of $\Delta$ that are in exactly one facet of $\Delta$. A face $F$ of $\Delta$ is {\bf interior} if $F$ is not in  $\partial\Delta$.

\section{Background on root polytopes}
\label{sec:root}

We follow the exposition of \cite[Section 4]{groth} in this section.  A root polytope of   type $A_{n}$ is  the convex hull of the origin and some points in $\Phi^+=\{e_i-e_j  \mid1\leq i<j \leq n+1\}$, the set positive roots of type $A_n$, where $e_i$ denotes the $i^{th}$ coordinate vector in $\mathbb{R}^{n+1}$. An important root polytope studied by Gelfand, Graev and Postnikov in  \cite{roothistory} is the full root polytope  $$\mathcal{P}(A_{n}^+)=\textrm{ConvHull}(0,  e_i-e_j \mid  1\leq i<j \leq n+1).$$ In this paper we restrict ourselves to a  class of root polytopes including $\mathcal{P}(A_{n}^+)$, which have subdivision algebras as defined in  \cite{root1}. We discuss subdivision algebras in relation to Grothendieck polynomials in Section \ref{sec:groth}.

Let $G$ be an acyclic   graph on the vertex set $[n+1]$.  Define $$\mathcal{V}_G=\{e_i-e_j \mid  (i, j) \in E(G), i<j\}, \mbox{ a set of vectors associated to $G$;}$$

 $${\rm cone}(G)=\langle \mathcal{V}_G \rangle :=\{\sum_{ e_i-e_j \in \mathcal{V}_G}c_{ij} (e_i-e_j) \mid  c_{ij}\geq 0\}, \mbox{ the cone associated to $G$; and } $$  
  $$\overline{\mathcal{V}}_G=\Phi^+ \cap {\rm cone}(G), \mbox{ all the positive roots of type $A_n$ contained in ${\rm cone}(G)$}.$$
    
The \textbf{root polytope} $\mathcal{P}(G)$ associated to the acyclic graph $G$ is  
 
 \begin{equation} \label{eq11} \mathcal{P}(G):=\textrm{ConvHull}(0, e_i-e_j \mid e_i-e_j \in \overline{\mathcal{V}}_G)\end{equation} The root polytope $\mathcal{P}(G)$ associated to a graph $G$ can also be defined as \begin{equation} \label{eq21} \mathcal{P}(G)=\mathcal{P}(A_n^+) \cap {\rm cone}(G).\end{equation}  
   
 Note that $\mathcal{P}(A_{n}^+)=\mathcal{P}(P_{n+1})$ for the  path graph  $P_{n+1}$ on the vertex set $[n+1]$.  
 
    The {\bf reduction rule for graphs:} Given   a graph $G_0$ on the vertex set $[n+1]$ and   $(i, j), (j, k) \in E(G_0)$ for some $i<j<k$, let   $G_1, G_2, G_3$ be graphs on the vertex set $[n+1]$ with edge sets
  \begin{eqnarray} \label{graphs}
E(G_1)&=&E(G_0)\backslash \{(j, k)\} \cup \{(i, k)\}, \nonumber \\
E(G_2)&=&E(G_0)\backslash \{(i, j)\} \cup \{(i, k)\},\nonumber \\ 
E(G_3)&=&E(G_0)\backslash \{(i, j), (j, k)\} \cup \{(i, k)\}.
\end{eqnarray}
    We say that $G_0$ \textbf{reduces} to $G_1, G_2, G_3$ under the reduction rules defined by equations (\ref{graphs}).
 

    \begin{lemma} \cite{root1} \label{reduction_lemma} \textbf{(Reduction Lemma for Root Polytopes)} 
Given   an acyclic  graph $G_0$ with $d$ edges let  $(i, j), (j, k) \in E(G_0)$ for some $i<j<k$ and $G_1, G_2, G_3$ as described by equations (\ref{graphs}).   Then  

$$\mathcal{P}(G_0)=\mathcal{P}(G_1) \cup \mathcal{P}(G_2)$$   where all polytopes  $\mathcal{P}(G_0), \mathcal{P}(G_1), \mathcal{P}(G_2)$ are   $d$-dimensional and    
$$\mathcal{P}(G_3)=\mathcal{P}(G_1) \cap \mathcal{P}(G_2)  \mbox{   is $(d-1)$-dimensional. } $$
\end{lemma}
      \medskip
      
       The Reduction Lemma says that performing a reduction on an acyclic graph $G_0$  is the same as dissecting the $d$-dimensional  polytope $\mathcal{P}(G_0)$ into two $d$-dimensional polytopes $\mathcal{P}(G_1)$ and $ \mathcal{P}(G_2)$, whose vertex  sets are  subsets of the vertex set of   $\mathcal{P}(G_0)$, whose interiors are disjoint, whose union is $\mathcal{P}(G_0)$, and whose intersection is a facet of both.   
       
       The following theorem in \cite{root1} describes a triangulation of the root polytope $\P(G)$ for any acyclic graph $G$. Its proof is based on the Reduction Lemma stated above.  We now define  the terminology used in the theorem. A graph $G$ on the vertex set $[n+1]$ is said to be \textbf{noncrossing} if there are no $1\leq i<j<k<l\leq n+1$ such that  $(i,k), (j, l)$ are edges of $G$. The graph $G$ is said to be \textbf{alternating} if at each vertex $v$ of $G$ all edges are either of the form $(v,i)$ for $v<i$ or of the form $(i,v)$ for $i<v$. Finally, the  directed transitive closure of the graph $G$ on the vertex set $[n+1]$, denoted by $\bar{G}$, is  $\bar{G}:=(V(G), E(G)\cup \{(i,j) \mid (i,j) \not \in E(G) \text{ and there exist } i<i_1< \ldots< i_k<j \text{ with } (i,i_1), (i_1, i_2), \ldots, (i_k, j) \in E(G)\}).$ 

\bt \label{thm:can} \cite{root1} Let $T_1, \ldots, T_k$ be the noncrossing alternating spanning trees of the directed transitive closure of the noncrossing acyclic graph $G$. Then $\P(T_1), \ldots, \P(T_k)$ are top dimensional simplices in a regular  triangulation of $\P(G)$.
\et

We note that there is a version of Theorem \ref{thm:can} in \cite{root1} that does not require the noncrossing condition on $G$; however, in the present paper we only invoke the above version which  has the advantage that it is  easier to state. 

We refer to the triangulation specified in Theorem \ref{thm:can} as the \textbf{canonical triangulation} of $\P(G)$. We remark that the polytopes  $\P(T_i)$, $i \in [k]$, in Theorem \ref{thm:can} are simplices, because the graphs $T_i$, $i\in [k]$, are alternating trees, and as such  they are their own transitive closure, with the roots corresponding to the edges of each of them linearly independent. Since each simplex $\P(T_i)$, $i \in [k]$, contains the vertex $0$, it follows that the canonical triangulation of $\P(G)$ also induces a triangulation of the vertex figure of $\P(G)$ at $0$, which we also call the canonical triangulation of the  vertex figure of $\P(G)$ at $0$. We sumarize some facts about the canonical triangulation in the following proposition.

 \begin{proposition}\label{prop:cantriang} Let $\C(G)$ denote the simplicial complex induced by the canonical triangulation of the vertex figure of $\P(G)$ at $0$. The vertices of $\C(G)$ are in bijection with edges $(i,j)$ in the directed transitive closure of $G$; the vertex of $\C(G)$ corresponding to $(i,j)$  is the intersection of the ray pointing to $e_i-e_j$ and the hyperplane by which we intersect $\P(G)$ to obtain the considered vertex figure.
 \end{proposition}



\section{Pipe dream complexes as triangulations of root polytopes}
\label{sec:constr}

In this section we give  geometric realizations of pipe dream complexes of permutations $\pi=1\pi'$, where $\pi'$ is dominant, in terms of triangulations of root polytopes. Indeed, we  construct  a  geometric realization of the subword complex that is the core of  $PD(1\pi')$.  To this end we start by defining a tree $T(\pi)$ for each permutation $\pi=1\pi'$,   $\pi'$ dominant.

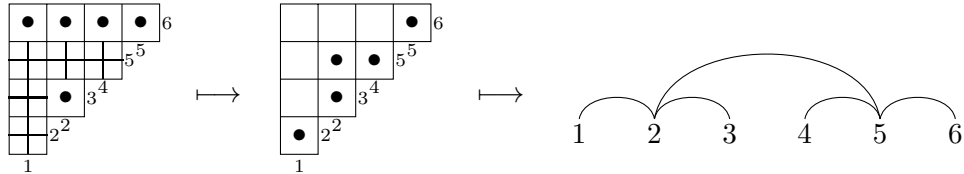
\begin{figure}[h]
    \centering
    \begin{subfigure}[h]{2.5cm}
        
          \begin{tikzpicture}
 
 \draw (0.5,2)--(0.5,0)--(0,0)--(0,2)--(2,2)--(2,1.5)--(0,1.5);
 \draw (0,0.5)--(1,0.5)--(1,2);
 \draw (0,1)--(1.5,1)--(1.5,2);

\draw (.25,1.75) node {$\bullet$};
\draw (.75,.75) node {$\bullet$};
\draw (.75,1.75) node {$\bullet$};
\draw (1.25,1.75) node {$\bullet$};
\draw (1.75,1.75) node {$\bullet$};
 
 \draw (0.26,0.15) node {$\textcross$};
 \draw (0.26,0.65) node {$\textcross$};
  \draw (0.26,1.15) node {$\textcross$};
   \draw (0.76,1.15) node {$\textcross$};
    \draw (1.26,1.15) node {$\textcross$};

\draw (0.25,-0.15) node {{\tiny 1}};
 \draw (.6,0.25) node {{\tiny 2}};
 \draw (.75,0.39) node {{\tiny 2}};
 \draw (1.1,0.75) node {{\tiny 3}};
 \draw (1.25,0.89) node {{\tiny 4}};
  \draw (1.6,1.25) node {{\tiny 5}};
 \draw (1.75,1.39) node {{\tiny 5}};
 \draw (2.1,1.75) node {{\tiny 6}};
 
\end{tikzpicture}

    \end{subfigure}
    \begin{subfigure}[h]{1cm}
        $\displaystyle{\lmp{}}$
        
            \end{subfigure}
        \begin{subfigure}[h]{2.5cm}
        
          \begin{tikzpicture}
 
 \draw (0.5,2)--(0.5,0)--(0,0)--(0,2)--(2,2)--(2,1.5)--(0,1.5);
 \draw (0,0.5)--(1,0.5)--(1,2);
 \draw (0,1)--(1.5,1)--(1.5,2);

\draw (.25,.25) node {$\bullet$};
\draw (.75,.75) node {$\bullet$};
\draw (.75,1.25) node {$\bullet$};
\draw (1.25,1.25) node {$\bullet$};
\draw (1.75,1.75) node {$\bullet$};

\draw (0.25,-0.15) node {{\tiny 1}};
 \draw (.6,0.25) node {{\tiny 2}};
 \draw (.75,0.39) node {{\tiny 2}};
 \draw (1.1,0.75) node {{\tiny 3}};
 \draw (1.25,0.89) node {{\tiny 4}};
  \draw (1.6,1.25) node {{\tiny 5}};
 \draw (1.75,1.39) node {{\tiny 5}};
 \draw (2.1,1.75) node {{\tiny 6}};
 
\end{tikzpicture}

    \end{subfigure}
    \begin{subfigure}[h]{1cm}
        $\displaystyle{\lmp{}}$
        
            \end{subfigure}
    \begin{subfigure}[h]{5.4cm}

        \begin{tikzpicture}
 
\path (1,0) edge [bend left=80] (2,0);
\path (2,0) edge [bend left=80] (3,0);
\path (2,0) edge [bend left=80] (5,0);
\path (4,0) edge [bend left=80] (5,0);
\path (5,0) edge [bend left=80] (6,0);

\draw[] (1,-0.2) node {1};
\draw[] (2,-0.2) node {2};
\draw[] (3,-0.2) node {3};
\draw[] (4,-0.2) node {4};
\draw[] (5,-0.2) node {5};
\draw[] (6,-0.2) node {6};
 
\end{tikzpicture}
        
    \end{subfigure}
\caption{Let $\pi=[15342]$, on the left we see the bottom reduced pipe dream for $\pi$ drawn inside $R(\pi)$ with dots instead of elbows, this gives the labeling of the boundary. We then drop the dots to the south to get the dots encoding $T(\pi)$, which is depicted on the right.}\label{fig:cvs}\end{figure}

Let $\pi=1\pi'$, where $\pi'$ is dominant. Denote by $\S(\pi)$ the subword complex that is the core$(\pi)$ coned over the vertex of $PD(\pi)$ corresponding to the entry $(1,1)$. Denote the region that is the union of $(1,1)$ and cr$(\pi)$ by $R(\pi)$. In order to determine the tree $T(\pi)$, we will label the southeast boundary with some numbers and we will place dots in some entries of $R(\pi)$, see Figure \ref{fig:cvs}. The boundary of the core region starting from the southwest (SW) corner of it to the northeast (NE) corner can be described as a series of east (E) and north (N) steps. Let $A$ be the set consisting of all the N steps together with some E steps. The step $\text{E}\in A$ if the bottom reduced pipe dream is bounded by E but not by the N step directly preceding E. As we traverse this lower boundary from the SW corner we write the numbers $1,\ldots, m$ in increasing fashion below the E steps and to the right of the N steps that belong to $A$. For the E steps that we did not assign a number, we consider their number to be the number assigned to the N step directly preceding them. 

We now describe how to place the dots in $R(\pi)$. Consider the bottom reduced pipe dream drawn inside $R(\pi)$ and with elbows replaced by dots. Drop these dots south. Define $T(\pi)$ to be the tree on $m$ vertices such that there is an edge between vertices $i<j$ if there is a dot in the entry in the column of the E step labeled $i$ and in the row of the N step labeled $j$. Let $t(\pi)$ be the number of edges of $T(\pi)$.

\noindent \textit{Remark.} We note that the above construction can be  simplified  in the special case of $\pi=1\pi'$ with $D(\pi')$ a partition with distinct parts. Indeed, in the above construction  the first string of E steps are in $A$. Furthermore, the E steps that can be seen as the boundary of the bottom reduced pipe dream B of $\pi$ and such that E bounds row $r$ of length $l_r$ of B and the row below row $r$ of length $l_{r+1}$ is at least two boxes shorter than row $r$ and moreover, step E is not $(l_{r+1}+1)^{st}$ from the left side are also in $A$. Also, we are placing dots in the rightmost boxes of the core region as well as in positions $(l_{r+1}+1)^{st}$ until $(l_r-1)^{th}$ in rows $r$ of the bottom reduced pipe dream  B that are longer than row $r+1$ by at least two boxes. In the case in which $\pi=1\pi'$ where $\pi'$ is dominant with its diagram having all parts distinct, then the decoration on the core diagram is much simpler. The first string of E steps consist of only one E step and this is also the only E step in $A$, so the number of vertices of $T(\pi)$ is 2 more than the size of the largest column of the diagram of $\pi$. The dots are placed on the rightmost boxes of cr$(\pi)$, see Figure \ref{dist}.

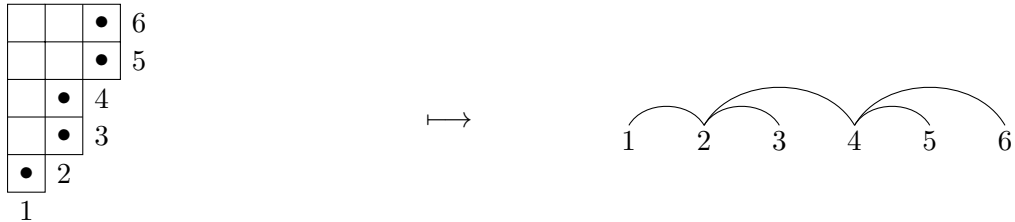
\begin{figure}[h]
    \centering
    \begin{subfigure}[h]{0.4\textwidth}
        \centering
        
        \begin{tikzpicture}
 
 \draw (0,1.5)--(1.5,1.5)--(1.5,2.5)--(0.5,2.5)--(0.5,0)--(0,0)--(0,2)--(1.5,2);
 \draw (0,0.5)--(1,0.5)--(1,2.5);
 \draw (0,1)--(1,1);
 \draw (0,2)--(0,2.5)--(0.5,2.5);
 
 \draw (0.25,0.25) node {$\bullet$};
\draw (0.75,0.75) node {$\bullet$};
\draw (0.75,1.25) node {$\bullet$};
\draw (1.25,1.75) node {$\bullet$};
\draw (1.25,2.25) node {$\bullet$};

\draw (0.25,-0.25) node {1};
\draw (0.75,0.25) node {2};
 \draw (1.25,0.75) node {3};
 \draw (1.25,1.25) node {4};
 \draw (1.75,1.75) node {5};
 \draw (1.75,2.25) node {6};

\end{tikzpicture}

    \end{subfigure}
    \hfill
    \begin{subfigure}[h]{0.1\textwidth}
        \centering   
        $\displaystyle{\lmp{}}$
        
            \end{subfigure}
    \hfill
    \begin{subfigure}[h]{0.4\textwidth}
        \centering

        \begin{tikzpicture}
 
\path (1,0) edge [bend left=60] (2,0);
\path (2,0) edge [bend left=60] (3,0);
\path (2,0) edge [bend left=60] (4,0);
\path (4,0) edge [bend left=60] (5,0);
\path (4,0) edge [bend left=60] (6,0);

\draw[] (1,-0.2) node {1};
\draw[] (2,-0.2) node {2};
\draw[] (3,-0.2) node {3};
\draw[] (4,-0.2) node {4};
\draw[] (5,-0.2) node {5};
\draw[] (6,-0.2) node {6};
 
\end{tikzpicture}

    \end{subfigure}
    \hfill
   
    \caption{Obtaining $T(\pi)$ from $R(\pi)$.}
    \label{dist}
\end{figure}

The vertices of $\S(\pi)$ are in bijection with configurations of one elbow and $|R(\pi)|-1$ crosses in $R(\pi)$, where $|R(\pi)|$ equals the number of entries in $R(\pi)$. Denote these vertices by $v_1, \ldots, v_k$. We define a map M from the vertices of the simplicial complex $\S(\pi)$ to the vertices of $\C(\pi):=\C(T(\pi))$. Recall that $\C(\pi)$ is the canonical triangulation of the vertex figure at $0$ of the root polytope $\P(T(\pi))$ and by Proposition \ref{prop:cantriang} the vertices of the triangulation are in bijection with edges $(i,j)$ in the directed transitive closure of $T(\pi)$.  The latter in turn are in bijection with the boxes of $R(\pi)$ by the map that takes a box to the edge $(i,j)$ if the E step below the box and in the boundary of $R(\pi)$  is labeled by $i$ and the N step to the right of the box and in the boundary of $R(\pi)$  is labeled by $j$. The map M is defined analogously as follows.
\begin{definition} Consider a vertex of $\S(\pi)$; this can be seen as a sole elbow tile in $R(\pi)$. M maps this vertex to the vertex of $\C(\pi)$ corresponding to $(i,j)$ if the box containing the elbow tile yields the edge $(i,j)$ in $T(\pi)$ (that is to the intersection of the ray pointing to $e_i-e_j$ and the hyperplane by which we intersect $\P(T(\pi))$ to obtain the considered vertex figure).
\end{definition}

\begin{figure}[h]
    \centering
    \begin{subfigure}[h]{0.4\textwidth}
        \centering
        
        \begin{tikzpicture}
 
 \draw (0,1.5)--(1.5,1.5)--(1.5,2.5)--(0.5,2.5)--(0.5,0)--(0,0)--(0,2)--(1.5,2);
 \draw (0,0.5)--(1,0.5)--(1,2.5);
 \draw (0,1)--(1,1);
 \draw (0,2)--(0,2.5)--(0.5,2.5);
 
 \draw (0.25,1.25) node {\color{red}$\bullet$};
\draw (0.75,1.75) node {\color{blue}$\bullet$};
\draw (0.75,1.25) node {\color{green}$\bullet$};

\draw (0.25,-0.25) node {1};
\draw (0.75,0.25) node {2};
 \draw (1.25,0.75) node {3};
 \draw (1.25,1.25) node {4};
 \draw (1.75,1.75) node {5};
 \draw (1.75,2.25) node {6};

\end{tikzpicture}

    \end{subfigure}
    \hfill
    \begin{subfigure}[h]{0.1\textwidth}
        \centering   
        $\displaystyle{\lmp{}}$
        
            \end{subfigure}
    \hfill
    \begin{subfigure}[h]{0.4\textwidth}
        \centering

        \begin{tikzpicture}
 
\path[color=red] (1,0) edge [bend left=60] (4,0);
\path[color=blue] (2,0) edge [bend left=60] (5,0);
\path[color=green] (2,0) edge [bend left=60] (4,0);

\draw[] (1,-0.2) node {1};
\draw[] (2,-0.2) node {2};
\draw[] (3,-0.2) node {3};
\draw[] (4,-0.2) node {4};
\draw[] (5,-0.2) node {5};
\draw[] (6,-0.2) node {6};
 
\end{tikzpicture}

    \end{subfigure}
    \hfill
   
    \caption{A few examples of the image of $M$. Each colored dot represents the vertex with an elbow at that position.}
    \label{dist}
\end{figure}
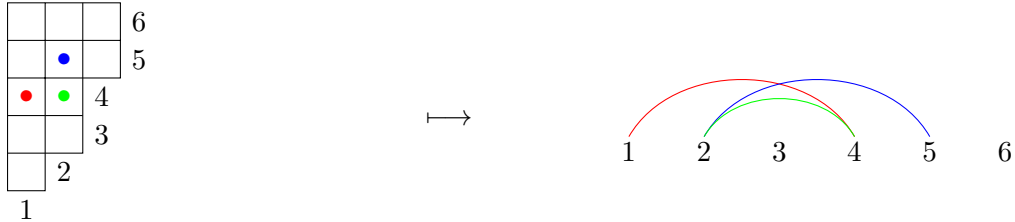

\bt \label{geom} The map M described above respects the simplicial complex structure of $\S(\pi)$ and $\C(\pi)$. In other words, $\C(\pi)$ is a  geometric realization of the subword complex $\S(\pi)$. \et

\proof Since both $\S(\pi)$ and $\C(\pi)$ are pure simplicial complexes of the same dimension by Lemma \ref{dim} it suffices to show that the map M is a bijection on the facets of $\S(\pi)$ and $\C(\pi)$. This is proven in Theorem \ref{thm:bij}.
\qed

\medskip

\noindent{\it Proof of Theorem \ref{thm:intro}.} This follows  using Theorem \ref{geom}. \qed

\medskip

Next we show that from Theorem \ref{geom} it follows that we can also realize  core$(\pi)$ geometrically, which is a subword complex as explained in Section \ref{sec:pipe} after Definition \ref{def:cone}. To this end we prove an auxiliary lemma first.

 \bl \label{int} Let $\pi=1\pi'$,  with  $\pi'$ dominant. If $\C(\pi)$ has an interior vertex, then it is the  unique vertex in $\C(\pi)$ on the ray between $0$ and $e_1-e_m$. Moreover, $\C(\pi)$ has an interior vertex if and only if $\pi=1n(n-1)\dots 2$.
 \el
 
 \proof  $\C(\pi)$ has an interior vertex if and only if the cone generated by $e_i-e_j$ for $(i,j)\in T(\pi)$ has an interior point $e_k-e_l$, where $(k,l)$ is in the directed transitive closure of $T(\pi)$. Since $e_i-e_j$ for $(i,j)\in T(\pi)$  are linearly independent, an interior point can be expressed as $\sum_{(i,j)\in T(\pi)} c_{ij} (e_i-e_j)$ with $c_{ij}>0$, $(i,j)\in T(\pi)$. If  $\mbox{ } T(\pi)=([m], \{(i,i+1) \mid i \in [m-1]\})$, then $e_1-e_m=\sum_{(i,j)\in T(\pi)} (e_i-e_j)$ is an interior point; moreover, $e_k-e_l$ for $1\leq k<l\leq m$ is an interior point if and only if $k=1$ and $l=m$. We have $T(\pi)=([m], \{(i,i+1) \mid i \in [m-1]\})$ exactly for  $\pi=1 m (m-1)\dots 2$. 
 For $T(\pi)\neq ([m], \{(i,i+1) \mid i \in [m-1]\})$, there is no $(k,l)$ in the directed transitive closure of $T(\pi)$ such that 
  $e_k-e_l=\sum_{(i,j)\in T(\pi)} c_{ij} (e_i-e_j)$  with $c_{ij}>0$.
\qed

\bt \label{thm:core} Let $\pi=1\pi'$,  with  $\pi'$ dominant. Let $v$ be the unique vertex in $\C(\pi)$ on the ray between $0$ and $e_1-e_m$.  For $\pi \neq 1 m (m-1)\dots 2$,  $v$ is in the boundary of $\C(\pi)$, and {\rm core}$(\pi)$
 is realized by the induced triangulation of the  vertex figure of $\C(\pi)$ at $v$. For $\pi= 1 m (m-1)\dots 2$,   {\rm core}$(\pi)$
 is realized by  the induced triangulation of the boundary of $\C(\pi)$.\et

\proof The vertex  $v$, which is the  unique vertex in $\C(\pi)$ on the ray between $0$ and $e_1-e_m$,  is the unique coning point of the geometric realization of $\S(\pi)$. If $v$ is in the boundary of  $\C(\pi)$,  which happens exactly when  $\pi \neq 1 m (m-1)\dots 2$ by Lemma \ref{int},  then the induced triangulation of  the vertex figure at   $v$ of $\C(\pi)$ is  a  geometric realization of 
 core$(\pi)$ (which is homeomorphic to a ball). For $\pi= 1 m (m-1)\dots 2$, the coning point $v$ lies in the interior of $\C(\pi)$, then since it is the only point in the interior of the canonical triangulation $\C(\pi)$ by Lemma \ref{int}, then the induced triangulation of the boundary of  $\C(\pi)$ is a  geometric realization of core$(\pi)$    (which is homeomorphic to a sphere). \qed

\bl  \label{dim} The core of the pipe dream complex of $\pi=1\pi'$, where the diagram of $\pi'$ is a partition, is of  dimension $t(\pi)-2$. The dimension of the root polytope $\P(T(\pi))$ is $t(\pi)$ and its vertex figure at $0$ is of dimension $t(\pi)-1$. In particular, both $\S(\pi)$ and $\C(\pi)$ are of dimension $t(\pi)-1$.  
\el
\proof
Since subword complexes are pure, then the dimension of the core of the pipe dream complex of $\pi$ equals the dimension of one of its facets. Consider the facet given by the bottom reduced pipe dream drawn inside the core region. The dimension of this facet equals one less than the number of elbows in the core and from the construction of $T(\pi)$ this equals $t(\pi)-2$. 
The dimension of  the root polytope $\P(T(\pi))$ is the number of edges in $T(\pi)$, which by definition is $t(\pi)$.
\qed

\medskip

The map $M$ can be easily extended to a map between  pipe dreams $P$ of $\pi$ drawn inside $R(\pi)$ and forests $F$  on $m$ vertices as follows. For each elbow tile in $P$ add the edge $(i,j)$ corresponding to the box of the elbow to $F$. Moreover, add the  edge $(1,m)$ to $F$. 

\bt \label{thm:bij} The reduced pipe dreams of $\pi=1\pi'$, where the diagram of $\pi'$ is a partition, are in bijection with the noncrossing alternating spanning trees of the directed transitive closure of $T(\pi)$ via the map $M$.
\et

We prove Theorem \ref{thm:bij} by induction on the number of columns in the diagram. We break it down in several lemmas.


\bl \label{lem:base} Take the permutation  $\pi=1\pi'$, where the diagram of $\pi'$ is $\lambda=(k)$. The reduced pipe dreams of $\pi$ are in bijection with the noncrossing alternating spanning trees of the directed transitive closure of $T(\pi)$ via the map M.
\el

\proof 

\begin{figure}[h]
    \centering
    \begin{subfigure}[h]{0.4\textwidth}
        \centering
        
        \begin{tikzpicture}
 
 \draw (0,0.5)--(0,0)--(0.5,0)--(0.5,2.5)--(1,2.5)--(1,0.5)--(0,0.5)--(0,2)--(1,2);
 \draw (0,1.5)--(1,1.5);
 \draw (0,1)--(1,1);
 
 \draw (0.25,0.25) node {$\bullet$};
\draw (0.75,0.75) node {$\bullet$};
\draw (0.75,1.25) node {$\bullet$};
\draw (0.75,1.75) node {$\bullet$};
\draw (0.75,2.25) node {$\bullet$};

\draw (0.25,-0.25) node {1};
\draw (0.75,0.25) node {2};
 \draw (1.25,0.75) node {3};
 \draw (1.25,1.25) node {4};
 \draw (1.25,1.75) node {5};
 \draw (1.25,2.25) node {6};

\end{tikzpicture}

    \end{subfigure}
    \hfill
    \begin{subfigure}[h]{0.1\textwidth}
        \centering   
        $\displaystyle{\lmp{M}}$
        
            \end{subfigure}
    \hfill
    \begin{subfigure}[h]{0.4\textwidth}
        \centering

        \begin{tikzpicture}
 
\path (1,0) edge [bend left=60] (2,0);
\path (2,0) edge [bend left=60] (3,0);
\path (2,0) edge [bend left=60] (4,0);
\path (2,0) edge [bend left=60] (5,0);
\path (2,0) edge [bend left=60] (6,0);

\draw[] (1,-0.2) node {1};
\draw[] (2,-0.2) node {2};
\draw[] (3,-0.2) node {3};
\draw[] (4,-0.2) node {4};
\draw[] (5,-0.2) node {5};
\draw[] (6,-0.2) node {6};
 
\end{tikzpicture}

    \end{subfigure}
    \hfill
    \caption{The core region and the tree $T(\pi)$ for $\lambda=(4)$} \label{fig:4}
\end{figure}
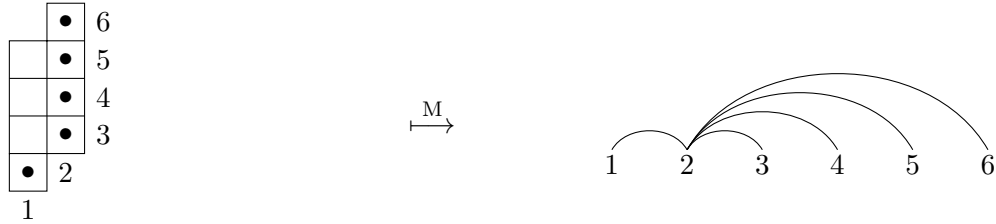

The edges of $T(\pi)$ for such a $\pi$ are $(1,2)$ and $(2,j)$ for $j=3,\ldots,k+2$ and thus for the transitive closure of $T(\pi)$ we add the edges $(1,j)$ with $j=3,\ldots,k+2$. See Figure \ref{fig:4}.

For $l=2,\ldots,k+2$ let $T_l$ be a tree on the vertex set $[k+2]$ consisting of the edges $(2,i)$ for $2<i\leq l$ and $(1,j)$ for $j\geq l$. Then $T_l$, $l=2,\ldots,k+2$, 
are all of the noncrossing alternating spanning trees of the directed transitive closure of $T(\pi)$. The map $M$ applied to the bottom reduced pipe dream of $\pi$ yields $T_{k+2}$. Moreover, after performing $0\leq i\leq k$ ladder moves (there is only one way to do this) on the bottom reduced pipe dream of $\pi$, we obtain a reduced pipe dream whose image under $M$ is $T_{k+2-i}$. By Theorem \ref{thm:ladder}  these are indeed all of the reduced pipe dreams of $\pi$ concluding the proof.
\qed


\bl \label{lem:lastcol} Given   $\pi=1\pi'$, where the diagram of $\pi'$ is a partition and that has more than one nonzero column, let its rightmost (shortest) column be of size $k$. Then in a reduced pipe dream of $\pi$ the only configurations of crosses and elbows that can occur in the rightmost column of {\rm cr}$(\pi)$ are, as read from above,  $l$ crosses and $k-l$ elbows, for $l=0,...,k$. 
\el

\proof This follows immediately from  Theorem \ref{thm:ladder}. \qed

\bl  \label{lem:gen} Let  $\pi=1\pi'$, where the diagram of $\pi'$ is a partition $\lambda=(\lambda_1,\ldots, \lambda_z)$ that has more than one nonzero column.  Consider all reduced pipe dreams of  $\pi$ where the configuration of crosses and elbows in the rightmost column of {\rm cr}$(\pi)$ is set to  have $l$ crosses and $k-l$ elbows for a fixed $0\leq l\leq k$. These are in bijection with  reduced pipe dreams of the permutation $1w_l$, where $w_l$ has diagram $(\lambda_1-(k-l), \lambda_2-(k-l), \ldots, \lambda_{z-1}-(k-l))$.
\el

\proof Since the bottom $k-l$ boxes of the rightmost column of cr$(\pi)$ are elbows, it can be seen using Theorem \ref{thm:ladder}   that the 
$k-l$ rows containing crosses one step to the south and one step to the west of these $k-l$ boxes can never move anywhere. 

Moreover, the fixed rows of crosses do not affect the ladder moves we can make on the remaining crosses. This allows us to get exactly the reduced pipe dreams for the permutation $1w$, where the diagram of $w$ is the diagram of $\pi$ after ignoring the fixed rows and shortest column, i.e., $w$ has diagram $(\lambda_1-(k-l), \lambda_2-(k-l), \ldots, \lambda_{z-1}-(k-l))$.
\qed

The following example illustrates the lemma above.

\begin{example} Let $1\pi=[164235]$ and suppose $l=1$, i.e., we are fixing one cross in position $(1,3)$ and elbow in position $(2,3)$, see Figure \ref{ex:fixelb}. The elbow in entry $(2,3)$ causes row 3 to consist of only crosses. Therefore the reduced pipe dreams for $[164235]$ with a cross in entry $(1,3)$ and elbow in $(2,3)$ correspond with the reduced pipe dreams for $[15234]$. The diagram of this latter permutation has fewer columns.

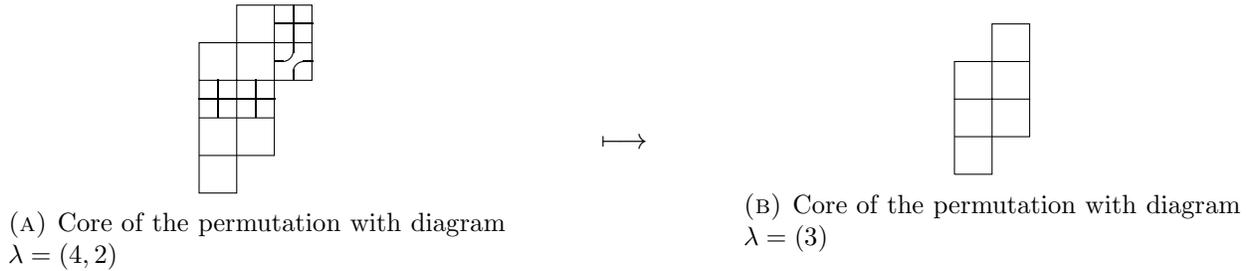
\begin{figure}[h]
    \centering
    \begin{subfigure}[h]{0.4\textwidth}
        \centering
        
        \begin{tikzpicture}
 
 \draw (0,1.5)--(1.5,1.5)--(1.5,2.5)--(0.5,2.5)--(0.5,0)--(0,0)--(0,2)--(1.5,2);
 \draw (0,0.5)--(1,0.5)--(1,2.5);
 \draw (0,1)--(1,1);
 
\draw (0.26,1.15) node {$\textcross$};
\draw (0.76,1.15) node {$\textcross$};
\draw (1.26,1.65) node {$\textelbow$};
\draw (1.26,2.15) node {$\textcross$};

\end{tikzpicture}

        \caption{Core of the permutation with diagram $\lambda=(4,2)$}
        \label{fig:y equals x}
    \end{subfigure}
    \hfill
    \begin{subfigure}[h]{0.1\textwidth}
        \centering   
        $\displaystyle{\lmp{}}$
        
            \end{subfigure}
    \hfill
    \begin{subfigure}[h]{0.4\textwidth}
        \centering    
        
        \begin{tikzpicture}
 
 \draw (0,0)--(.5,0)--(0.5,2)--(1,2)--(1,1)--(0,1)--(0,0);
 \draw (0,1)--(0,1.5)--(1,1.5);
\draw (0,.5)--(1,.5)--(1,1);

\end{tikzpicture}      
        
        \caption{Core of the permutation with diagram $\lambda=(3)$}
    \end{subfigure}
    \hfill
    \caption{New core after applying the reduction of Lemma \ref{lem:gen} for $l=1$ to the core on the left.} \label{ex:fixelb} 
\end{figure}
\end{example}
  
\medskip

\bl \label{lem:S} Given $\pi=1\pi'$, $\pi'$ dominant, where the length of the shortest column of the diagram of $\pi'$ is $k$,  the set $S$ of noncrossing alternating spanning trees of the directed transitive closure of $T(\pi)$ is a disjoint union $S=S_0\sqcup \cdots \sqcup S_{k}$, where 
	\begin{align}S_l= &\{T\in S : (m-k,m-j)\notin E(T) \text{ for } j=0,\ldots,l-1\}\\
	& \cup \{T\in S : (m-k,m-j)\in E(T) \text{ for } j=l,\ldots,k-1\},\notag\end{align}
for $0\leq l\leq k$, where $m$ is the number of vertices of $T(\pi)$.
\el

Note that $m-k$ is the label on the bottom of the last column of the core of $\pi$. Thus, $S_l$ consists of the noncrossing alternating spanning trees of the directed transitive closure of $T(\pi)$ that do not contain the edges corresponding to the top $l$ crosses in last column of the core of $\pi$ and contain the edges corresponding to the bottom $k-l$ elbows in last column of the core.

\noindent {\it Proof of Lemma \ref{lem:S}.} This follows immediately from the definition of $T(\pi)$.
\qed

\bl\label{lem:fixing} Let $\pi=1\pi'$, $\pi'$ dominant of shape $(\lambda_1, \ldots, \lambda_z)$ and let  $1w_l$ be the permutation where $w_l$ has diagram $(\lambda_1-(k-l), \lambda_2-(k-l), \ldots, \lambda_{z-1}-(k-l))$. Use Lemma \ref{lem:gen} to draw the core region of $1w_l$ inside the core region of $\pi$. Then all the edges corresponding to the entries outside the core region of $1w_l$ in a tree $T\in S_l$ are forced by the last column.
\el

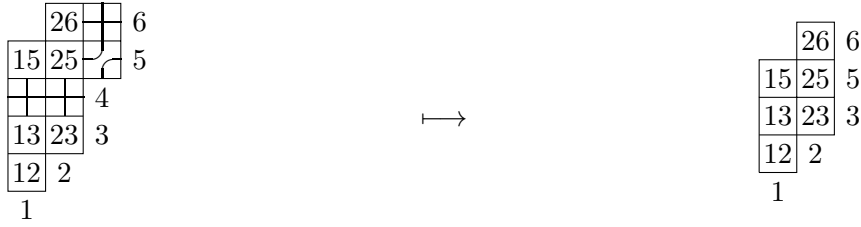
\begin{figure}[h]
    \centering
    \begin{subfigure}[h]{0.4\textwidth}
        \centering
        
        \begin{tikzpicture}
 
 \draw (0,1.5)--(1.5,1.5)--(1.5,2.5)--(0.5,2.5)--(0.5,0)--(0,0)--(0,2)--(1.5,2);
 \draw (0,0.5)--(1,0.5)--(1,2.5);
 \draw (0,1)--(1,1);
 
\draw (0.26,1.15) node {$\textcross$};
\draw (0.76,1.15) node {$\textcross$};
\draw (1.26,1.65) node {$\textelbow$};
\draw (1.26,2.15) node {$\textcross$};

\draw (0.25,.25) node {12};
\draw (0.25,.75) node {13};
\draw (0.75,.75) node {23};
\draw (0.25,1.75) node {15};
\draw (0.75,1.75) node {25};
\draw (0.75,2.25) node {26};

\draw (0.25,-.25) node {1};
\draw (0.75,.25) node {2};
\draw (1.25,.75) node {3};
\draw (1.25,1.25) node {4};
\draw (1.75,1.75) node {5};
\draw (1.75,2.25) node {6};

\end{tikzpicture}

        \label{fig:y equals x}
    \end{subfigure}
    \hfill
    \begin{subfigure}[h]{0.1\textwidth}
        \centering   
        $\displaystyle{\lmp{}}$
        
            \end{subfigure}
    \hfill
    \begin{subfigure}[h]{0.4\textwidth}
        \centering    
        
        \begin{tikzpicture}

 \draw (0,0)--(.5,0)--(0.5,2)--(1,2)--(1,1)--(0,1)--(0,0);
 \draw (0,1)--(0,1.5)--(1,1.5);
\draw (0,.5)--(1,.5)--(1,1);
 
 \draw (0.25,-.25) node {1};
\draw (0.75,.25) node {2};
\draw (1.25,.75) node {3};
\draw (1.25,1.25) node {5};
\draw (1.25,1.75) node {6};

\draw (0.25,.25) node {12};
\draw (0.25,.75) node {13};
\draw (0.75,.75) node {23};
\draw (0.25,1.25) node {15};
\draw (0.75,1.25) node {25};
\draw (0.75,1.75) node {26};

\end{tikzpicture}      
        
    \end{subfigure}
    \hfill
    \caption{Edge labeling for the core of $1w_1$ coming from the core of $[164235]$ on the left.} \label{fig:labelling} 
\end{figure}

\proof If $\lambda_z<\lambda_{z-1}$, then the edges inside cr$(\pi)$ and outside of cr$(1w_l)$ are precisely those in the last column and in the $(k-l)$ rows one step to the south of the $k-l$ boxes fixed to be elbows. The boxes in these rows are crosses and thus we conclude that in this case all the edges corresponding to boxes outside cr$(1w_l)$ are indeed fixed after fixing the last column.
If  $\lambda_z=\lambda_{z-1}=\cdots=\lambda_{z-j}<\lambda_{z-j-1}$, then aside from the boxes outside of cr$(1w_l)$ described in the previous sentence, we also have the $j$ boxes to the left of the top most elbow on the last column. We will show that a tree $T\in S_l$ must contain the edge corresponding to these boxes. Let $T\in S_l$, $u$ be the E step below the leftmost of these boxes and $v$ be the N step to the right of these boxes. Since $T$ is an alternating spanning tree, then $v$ must be adjacent either to $u$ or to a vertex before $u$. Similarly, $u$ must be adjacent either to $v$ or to a vertex after $v$. The only way noncrossing is preserved is if $(u,v)$ is an edge of $T$. We continue in this fashion by looking at the second leftmost box and taking the E step below it and again prove that the edge corresponding to that box is in $T$.
\qed

\bl  \label{lem:l} The set $S_l$, $0\leq l\leq k$, as in Lemma \ref{lem:S} is in bijection with reduced pipe dreams of the permutation $1w_l$, where $w_l$ has diagram $(\lambda_1-(k-l), \lambda_2-(k-l), \ldots, \lambda_{z-1}-(k-l))$, via the map $M$.
\el

We prove this lemma and Theorem \ref{thm:bij} together by using induction on the number of columns of the diagram of $\pi'$.

\medskip

\noindent \textit{Proof of Theorem \ref{thm:bij} and Lemma \ref{lem:l}.} We use induction on the number of columns of the diagram of $\pi'$. The base case for Theorem \ref{thm:bij} where this diagram contains one  column is proved in Lemma \ref{lem:base}. Notice that in the proof of this lemma the base case for Lemma \ref{lem:l} is also proven.

By Lemma \ref{lem:S} the noncrossing alternating spanning trees of the directed transitive closure of $T(\pi)$ can be broken down into the sets $S_l$, $l=0,1,\ldots, k$. Consider the permutation $1w_l$ where $w_l$ has diagram $(\lambda_1-(k-l), \lambda_2-(k-l), \ldots, \lambda_{z-1}-(k-l))$. By inductive hypothesis, we know that $1w_l$ satisfies Theorem \ref{thm:bij}, i.e., its reduced pipe dreams are in bijection with the noncrossing alternating spanning trees of the directed transitive closure of $T(1w_l)$ via the map M. By Lemma \ref{lem:fixing}   we have that the noncrossing alternating spanning trees of the directed transitive closure of $T(1w_l)$ yield the set $S_l$. 
Finally, by Lemma \ref{lem:gen} we know that the reduced pipe dreams of the permutations $1w_l$, as $l=0,1,\ldots, k$,  are in bijection with the reduced pipe dreams of $\pi$, concluding the proof.
\qed

\medskip

\section{Reduced forms in the subdivision algebra and Grothendieck polynomials}
\label{sec:groth}
In this section we show that Grothendieck polynomials of permutations $\pi=1\pi'$, $\pi'$ dominant, are special cases of reduced forms in the subdivision algebra of root polytopes. To this end we start by defining the notions appearing in the previous sentence.

The \textbf{subdivision algebra of root polytopes} $\t$ is a commutative algebra generated by the variables $x_{ij}$, $1\leq i<j\leq n$,  over $\mathbb{Q}[\beta]$, subject to the relations $x_{ij} x_{jk}=x_{ik}(x_{ij}+x_{jk}+\beta)$, for $1\leq i<j<k\leq n$.  This algebra is called the subdivision algebra,  because its relations can be seen geometrically  as subdividing  root polytopes via Lemma \ref{reduction_lemma}. The subdivision algebra has been used extensively for subdividing  root (and flow) polytopes  in \cite{prod, mm, h-poly1, h-poly2, groth, root1, root2}.

 A  \textbf{reduced form} of the monomial   in the algebra $\t$  is a polynomial   obtained by successively substituting $x_{ik}(x_{ij}+x_{jk}+\beta)$ in place of an occurrence of $x_{ij} x_{jk}$ for some $i<j<k$ until no further reduction is possible. Note that the reduced forms are not necessarily unique.
  
A possible sequence of reductions in algebra $\t$ yielding a reduced form of $x_{12}x_{23}x_{34}$ is given by

\begin{eqnarray} \label{ex1}
 x_{12} \mbox {\boldmath$  x_{23}x_{34}$} & \rightarrow & \mbox{\boldmath$x_{12}$}x_{24}\mbox{\boldmath$x_{23}$}+\mbox{\boldmath$x_{12}$}x_{34}\mbox {\boldmath$x_{24}$}+\beta \mbox {\boldmath$x_{12}x_{24}$} \nonumber \\
& \rightarrow& \mbox {\boldmath$ x_{24}$} x_{13}\mbox {\boldmath$x_{12}$}+x_{24}x_{23}x_{13}+  \beta x_{24}x_{13}+x_{34}x_{14}x_{12}+x_{34}x_{24}x_{14} \nonumber \\
& &+\beta x_{34}x_{14}+\beta x_{14}x_{12}+\beta x_{24}x_{14}+\beta^2 x_{14} \nonumber \\
& \rightarrow  &x_{13}x_{14}x_{12}+x_{13}x_{24}x_{14}+\beta x_{13}x_{14}+x_{24}x_{23}x_{13}+\beta x_{24}x_{13}\nonumber \\
& & +x_{34}x_{14}x_{12}+x_{34}x_{24}x_{14}+\beta x_{34}x_{14}+\beta x_{14}x_{12}+\beta x_{24}x_{14}\nonumber \\
& &+\beta^2 x_{14}
\end{eqnarray}

\noindent where the pair of variables on which the reductions are performed is in boldface. The reductions are performed on each monomial separately.

Given a noncrossing tree $T$ on the vertex set $[n]$, let $m[T]:=\prod_{(i,j)\in T} x_{ij}$. The \textbf{canonical reduced form} ${\rm Crf}_T(x_{ij} \mid 1\leq i<j\leq n)$ of $m[T]$ is the reduced form obtained by performing reductions on the tree $T$ from front to back (or back to front) on the topmost edges always. This can of course be translated into an algebraic context as follows. For $x_{ij}=t_i$ we denote by  ${\rm Crf}_T(t_{1}, \ldots, t_{n-1})={\rm Crf}_T(x_{ij} \mid 1\leq i<j\leq n)$.  While the reduced form of a monomial in the subdivision algebra is not necessarily unique, once we set $x_{ij}=t_i$ it becomes unique. This is the statement of the next theorem which we prove in Section \ref{sec:unique}.

\bt \label{thm:unique} Given a noncrossing tree $T$ on the vertex set $[n]$, let $R_T(x_{ij} \mid 1\leq i<j\leq n)$ be an arbitrary reduced form of $m[T]$. Let $R_T(t_{1}, \ldots, t_{n-1})$ be the reduced form  $R_T(x_{ij} \mid 1\leq i<j\leq n)$ when we let $x_{ij}=t_i$. Then, 

\be R_T(t_{1}, \ldots, t_{n-1})={\rm Crf}_T(t_{1}, \ldots, t_{n-1}).\ee
\et

We will use the notation $\tilde{R}_T({\bf t})$ when instead of setting $x_{ij}=t_i$, we do the following. Let $i_1< \ldots< i_v$ be the vertices of $T$ that have outgoing edges. Therefore, the only $x_{ij}$'s appearing in a reduced form must have $i \in \{i_1, \ldots, i_v\}$.  The reduced form $\tilde{R}_T({\bf t})$ is then obtained from $R_T(x_{ij} \mid 1\leq i<j\leq n)$ by setting $x_{i_k, j}=t_k$ for $k \in [v]$ and all $j\in [n]$. 

\medskip

The following theorem provides a combinatorial way of thinking about double  Grothendieck polynomials.

\bt \label{allen} \cite{subword, fom-kir} The \textbf{double  Grothendieck polynomial $\mathfrak{G}_w({\bf x, y})$} for $w \in S_n$, where ${\bf x}=(x_1, \ldots, x_{n-1})$ and ${\bf y}=(y_1, \ldots, y_{n-1})$ can be written as

\begin{equation} \mathfrak{G}_w({\bf x, y})=\sum_{P\in {\rm Pipes}(w)}(-1)^{codim_{PD(w)}F(P)} wt_{x,y}(P),
\label{eq:groth}
\end{equation}

\noindent where  ${\rm Pipes}(w)$ is the set of all pipe dreams of $w$ (both reduced and nonreduced), $F(P)$ is the interior face in $PD(w)$ labeled by the pipe dream $P$,  $codim_{PD(w)} F(P)$ denotes the codimension of $F(P)$ in $PD(w)$ and  $wt_{x,y}(P)=\prod_{(i,j) \in {\rm cross}(P)} (x_i-y_j+x_iy_j)$, with ${\rm cross}(P)$ being the set of positions where $P$ has a cross. \et

Note that in the product $\prod_{(i,j) \in {\rm cross}(P)} (x_i-y_j+x_iy_j)$ appearing in the statement of Theorem \ref{allen} we are assuming a certain labeling of rows and columns. Conventionally, rows are labeled increasingly from top to bottom and columns are labeled  increasingly from left to right.  Also recall that the lowest degree terms of $\mathfrak{G}_w({\bf x, y})$ give the Schubert polynomial $\mathfrak{S}_w({\bf x, y})$. Except in Theorem \ref{thm:gen}, we will be working with the single Grothendieck  polynomial $$\mathfrak{G}_w({\bf x}):=\mathfrak{G}_w({\bf x, 0}).$$ In other words, for single Grothendieck polynomials we use the  weight $wt_{x}(P)=\prod_{(i,j) \in {\rm cross}(P)} x_i$ instead of $wt_{x,y}(P)=\prod_{(i,j) \in {\rm cross}(P)} (x_i-y_j+x_iy_j)$ in equation \eqref{eq:groth}.


In the spirit of Theorem \ref{allen}, we use the following definition for the {\bf $\beta$-Grothendieck polynomial}:

\begin{equation} \mathfrak{G}^{\beta}_w({\bf x}):=\sum_{P\in {\rm Pipes}(w)}\beta^{codim_{PD(w)}F(P)} wt_{x}(P).
\label{eq:bgroth}
\end{equation}

If we set $\beta=0$ in \eqref{eq:bgroth}, then we recover the single Schubert polynomial $\mathfrak{S}_w({\bf x})$. Note that if  in \eqref{eq:bgroth} we assume that $\beta$ has degree $-1$, while all other variables are of degre $1$, then the powers of $\beta$'s simply  make the polynomial $ \mathfrak{G}_w^{\beta}({\bf x})$ homogeneous. We chose this definition of $\beta$-Grothendieck polynomials, as it will be the most convenient notationwise for our purposes.

\bt  \label{thm:t} Given $\pi=1\pi'$, $\pi'$ dominant, we have that for any reduced form of $m[T(\pi)]$

\be \tilde{R}_{T(\pi)}({\bf t})=\big(\prod_{i=1}^{n-1} t_i^{g_i}\big) \mathfrak{G}_{{\pi}^{-1}}^{\beta}(t_1^{-1}, \ldots, t_{n-1}^{-1}),\ee
where $g_i$ is the number of boxes in the $i$th column from the left in $R(\pi)$.
\et

A special case of Theorem \ref{thm:t} for $\pi=1n(n-1)\cdots 2$ appears in \cite{k2014}.

  We now relate the canonical reduced form to the double Grothendieck polynomial. For $\beta=-1$ denote the canonical reduced form by ${\rm Crf}^{\beta=-1}_T(x_{ij} \mid 1\leq i<j\leq n)$.  Before we can state and prove Theorem \ref{thm:gen} we need to define a map $\phi$ from the labels $(i,j)$ that the boxes in the region $R(\pi)$ inherit from the labeling of its boundary (as described in  Figure \ref{fig:cvs})  to the conventional labeling where rows are labeled increasingly from top to bottom and columns are labeled  increasingly from left to right.  We call the former labeling the tree labeling and when unclear which labeling we are talking of we put a $T$ index on it: $(i,j)_T$. The map $\phi$ simply takes the tree label $(i,j)$ to the conventional label $(\phi_{ij}(i), \phi_{ij}(j))$ of the corresponding box. In the example of  Figure \ref{fig:cvs} we have that $\phi((1,6))=(1,1), \phi((2,3))=(3,2), \phi((5,6))=(1,4), \phi((4,5))=(2,3)$, and so forth.

\bt  \label{thm:gen} Given $\pi=1\pi'$, $\pi'$ dominant, we have that

\be  \label{eq:crf}   {\rm Crf}^{\beta=-1}_{T(\pi)}(x_{ij}=\frac{1}{x_{\phi_{ij}(i)}-y_{\phi_{ij}(j)}+x_{\phi_{ij}(i)}y_{\phi_{ij}(j)}} \mid 1\leq i<j\leq n)=\ee \be \nonumber =(\prod_{(i,j)_T \in R(\pi)} \frac{1}{x_{\phi_{ij}(i)}-y_{\phi_{ij}(j)}+x_{\phi_{ij}(i)}y_{\phi_{ij}(j)}}) \mathfrak{G}_{\pi}({\bf x, y}).\ee
\et

\proof By definition we have that  \be {\rm Crf}_{T(\pi)}(x_{ij} \mid 1\leq i<j\leq n)=\sum_{G \in {\mathcal{L}}(T(\pi))}\beta^{codim_{\P(T(\pi))}\P(G)} wt(G),\ee where $wt(G)=\prod_{(i,j)\in G}x_{ij}$, ${\mathcal{L}}(T(\pi))$ denotes the set of graphs corresponding to the terms of the reduced form of $m[T(\pi)]$, and $\P(G)$ denotes the simplex in the canonical triangulation of $\P(T(\pi))$ corresponding to $G$. Together with Theorem \ref{geom} using the map M and Theorem \ref{allen}, we obtain \eqref{eq:crf}.
\qed

\section{Volumes and Ehrhart series of root polytopes}
\label{sec:cor}

In this section we state the two immediate corollaries regarding volumes and Ehrhart series of root polytopes following from Theorem \ref{thm:intro}. Recall that the  normalized volume of a $d$-dimensional polytope
$\P$ is $d!$ times its usual volume, which is always integral for lattice polytopes.

\bt \label{thm:vol} Let $\pi=1\pi'$, where $\pi'$ is a dominant permutation. Then the normalized volume of $\P(T(\pi))$ is equal to the number of reduced pipe dreams of $\pi$. This can be written as
\be {\rm vol}(\P(T(\pi)))=\mathfrak{G}^{\beta=0}_{\pi} ({\bf 1}).\ee
\et

Recall that for a  polytope $\mathcal{P}\subset \mathbb{R}^{N}$, the {$t^{th}$ dilate} of $\mathcal{P}$ is $\displaystyle t \mathcal{P}=\{(tx_1, \ldots, tx_{N}) \mid  (x_1, \ldots, x_{N}) \in \mathcal{P}\}.$ The number of lattice points of $t\mathcal{P}$, where $t$ is a nonnegative integer and $\mathcal{P}$ is a convex polytope, is given by the {\bf Ehrhart function} $i({\mathcal{P}}, t)$. If $\mathcal{P}$ has integral vertices then $i({\mathcal{P}}, t)$ is a  polynomial.

In order to state the Ehrhart series of root polytopes we need the following lemma, which follows from the well-known relationship of $f$- and $h$-vectors. We note that we take $h(\C, x)=\sum_{i=0}^d h_{i}x^{i}$ to be the  $h$-polynomial of a $(d-1)$-dimensional simplicial complex $\C$.

\begin{lemma} \cite{Stcom} \label{pure} Let $\C$ be a $(d-1)$-dimensional pure  simplicial complex homeomorphic to a ball and $f_i^\circ$ be the number of interior faces of $\C$ of dimension $i$. Then 
\begin{equation} \label{h} h(\C, \beta+1)=\sum_{i=0}^{d-1} f_i^\circ \beta^{d-1-i}
\end{equation}
\end{lemma}

We also use:
\begin{lemma} \label{lem:groth} \cite{groth} For any permutation $\pi$ the following holds: 
$$\mathfrak{G}_{\pi}^{\beta-1}({\bf 1})=h(PD(\pi), \beta).$$
\end{lemma}

\begin{theorem} \label{thm-ehr} Let $\pi=1\pi'$, where $\pi'$ is a dominant permutation.
Then
\be \label{h-ehrhart} \mathfrak{G}_{\pi}^{\beta-1}({\bf 1})=\sum_{m\geq 0}(i(\P(T(\pi)), m)\beta^m)(1-\beta)^{\dim(\P(T(\pi)))+1}.\ee
\end{theorem}

\proof Since the canonical  triangulation $\C$ of $\P(T(\pi))$   is unimodular, we have \cite{Stcom} 
\be \label{h-ehrhart1} h(\C, \beta)=\sum_{m\geq 0}(i(\P(T(\pi)), m)\beta^m)(1-\beta)^{\dim(\P(T(\pi)))+1}.\ee By Theorem \ref{thm:intro} and Lemma \ref{pure} we get that  $h(\C, \beta)=h(PD(\pi), \beta)$. Together with Lemma \ref{lem:groth} this concludes the proof.
\qed

\section{Uniqueness of $t$-reduced forms}
\label{sec:unique}

The aim of this section is to prove Theorem \ref{thm:unique}, which states that when we let $x_{ij}=t_i$ for all $i$, then the reduced form becomes unique. For clarity we call the reduced forms with the substitution $x_{ij}=t_i$, the \textbf{t-reduced forms}.

In order to prove Theorem \ref{thm:unique} we recall several definitions and results from \cite{root1}.

A reduction on the edges $(i,j),(j,k)$ of a noncrossing graph $G$  is \textbf{noncrossing} if the graphs resulting from the reduction are also noncrossing. Analogously we can define noncrossing reductions on $m[G]$.
 
 \bt \label{thm:noncross}\cite{root1} Let $T$ be  a noncrossing tree on the vertex set $[n]$. Performing noncrossing reductions on $m[T]$, regardless of order,  we obtain a unique reduced form $R^{noncross}_T(x_{ij} \mid 1\leq i<j\leq n)$  for $m[T]$.
 \et

\begin{figure}[htbp] 
\begin{center}         \begin{tikzpicture}
 
\path (1,0) edge [bend left=70] (5,0);
\path (2,0) edge [bend left=60] (5,0);
\path (3,0) edge [bend left=60] (4,0);
\path (4,0) edge [bend left=60] (5,0);
\path (5,0) edge [bend left=60] (6,0);
\path (5,0) edge [bend left=60] (8,0);
\path (6,0) edge [bend left=60] (7,0);

\draw[] (1,-0.2) node {1};
\draw[] (2,-0.2) node {2};
\draw[] (3,-0.2) node {3};
\draw[] (4,-0.2) node {4};
\draw[] (5,-0.2) node {5};
\draw[] (6,-0.2) node {6};
\draw[] (7,-0.2) node {7};
\draw[] (8,-0.2) node {8};
 
\end{tikzpicture}
\caption{The edge sets of the pseudo-components in the graph depicted are $\{(1, 5), (5, 8)\}, \{(2, 5)\}, \{(3, 4), (4, 5)\}, \{(5, 6), (6, 7)\}$. The pseudo-component with edge set  $\{(1, 5), (5, 8)\}$ is a both a left and right pseudo-component, while the  pseudo-components with edge sets     $ \{(2, 5)\}, \{(3, 4), (4, 5)\}$     are  left  pseudo-components and the  pseudo-component with edge set  $\{(5, 6), (6, 7)\}$ is a  right pseudo-component.} 
\label{fig:pseudo}
\end{center} 
\end{figure}
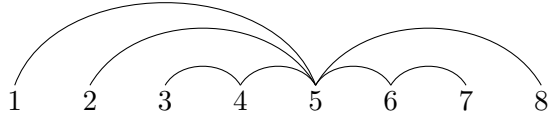

Consider a noncrossing tree $T$ on $[n]$. We define the \textbf{pseudo-components} of $T$ inductively. The unique simple path $P$ from $1$ to $n$ is a pseudo-component of $T$. The graph $T\backslash P$ is an edge-disjoint union of   trees $T_1, \ldots, T_k$, such that   if $v$ is a vertex of $P$ and $v \in T_l$, $l \in [k]$, then $v$ is either the minimal or maximal vertex of $T_l$ . Furthermore, there are no $k-1$ trees whose edge-disjoint  union is $T\backslash P$ and that satisfy  all the requirements stated above. The set of pseudo-components of $T$, denoted by $ps(T)$ is $ps(T)=\{P\} \cup ps(T_1)\cup \cdots \cup ps(T_k)$. A pseudo-component $P'$ is said to be on $[i, j]$, $i<j$  if it is a path with endpoints $i$ and  $j$. A pseudo-component $P' $ on $[i, j]$ is said to be a \textbf{left pseudo-component} of $T$  if there are no edges $(s, i) \in E(T)$ with $s<i$   and a \textbf{right pseudo-component}  if  if there are no edges $(j, s) \in E(T)$ with $j<s$.     See Figure \ref{fig:pseudo} for an example.
   
\begin{theorem} \label{gen:non} \cite{root1} Let $T$ be a noncrossing tree. Then $R^{noncross}_T(x_{ij} \mid 1\leq i<j\leq n)$ is the sum of the monomials corresponding to the following graphs weighted with powers of  $\beta$ (of degree $1$) to obtain a homogeneous polynomial. The graphs are: all 
 noncrossing alternating spanning  forests of the directed transitive closure of $T$  on the vertex set $[n]$ containing edge $(1, n)$ and at least one edge of the form $(i_1, j)$ with $i_1\leq i$ for each right  pseudo-component of $T$ on $[i, j]$ and  at least one edge of the form $(i, j_1)$ with $j\leq j_1$ for each left pseudo-component of $T$ on $[i, j]$. 
 \et
 
 We note that in the above we assume that the vertices of our graphs are drawn on a line from left to right in increasing order, $1, 2, \ldots, n.$ This condition is of course not an essential condition for the above theorems, and if we rearrange the order of the vertices of our graphs,  then we can reinterpret the above results accordingly.

  Consider the noncrossing tree $T$ on the vertex set $[n]$ with vertices drawn from left to right in increasing order $1, 2,\ldots, n$. Let $(k,l), (l,m)$ be a pair of nonalternating edges in $T$. If  the reduction performed on $(k,l), (l,m)$ is noncrossing, then we set $T_{klm}=T$ with $T_{klm}$ drawn identically to $T$. If  the reduction performed on $(k,l), (l,m)$ is not noncrossing, then let $C_1=(V_1, E_1)$ and $C_2=(V_2, E_2)$ be the connected components containing the vertices $k$ and $m$, respectively,  in the graph $T-\{(k,l), (l,m)\}=([n], E(T)\backslash \{(k,l), (l,m)\})$. Then we define $T_{klm}=T$ to be drawn with its vertices arranged from left to right in the following order: $v_1^1, \ldots, v_1^p, w_1, \ldots, w_q, v_2^1, \ldots, v_2^r$, where $V_1=\{v_1^1< \cdots< v_1^p\}$, $V_2=\{v_2^1< \cdots< v_2^r\}$, $ [n]\backslash (V_1 \cup V_2)=\{w_1< \cdots< w_q,\}.$ The tree  $T_{klm}$ is then a  noncrossing tree. See Figure \ref{fig:tklm}.  
    
  \begin{figure}[h]
    \centering
    \begin{subfigure}[h]{0.5\textwidth}
        \centering
        
\begin{tikzpicture}
 
\path (1,0) edge [bend left=70] (4,0);
\path (2,0) edge [bend left=60] (3,0);
\path (2,0) edge [bend left=60] (4,0);
\path (4,0) edge [bend left=60] (5,0);
\path (5,0) edge [bend left=60] (6,0);
\path (5,0) edge [bend left=60] (7,0);
\path (4,0) edge [bend left=60] (8,0);

\draw[] (1,-0.2) node {1};
\draw[] (2,-0.2) node {2};
\draw[] (3,-0.2) node {3};
\draw[] (4,-0.2) node {4};
\draw[] (5,-0.2) node {5};
\draw[] (6,-0.2) node {6};
\draw[] (7,-0.2) node {7};
\draw[] (8,-0.2) node {8};

\draw[] (4.5,-0.7) node {$T$};

\end{tikzpicture}     
        
    \end{subfigure}
    
    \vspace{0.5cm}
    
    \begin{subfigure}[h]{0.5\textwidth}
        \centering    
        
\begin{tikzpicture}
 
\path (1,0) edge [bend left=70] (2,0);
\path (1,0) edge [bend left=60] (4,0);
\path (3,0) edge [bend left=60] (4,0);
\path (4,0) edge [bend left=60] (5,0);
\path (4,0) edge [bend left=60] (6,0);
\path (6,0) edge [bend left=60] (7,0);
\path (6,0) edge [bend left=60] (8,0);

\draw[] (1,-0.2) node {2};
\draw[] (2,-0.2) node {3};
\draw[] (3,-0.2) node {1};
\draw[] (4,-0.2) node {4};
\draw[] (5,-0.2) node {8};
\draw[] (6,-0.2) node {5};
\draw[] (7,-0.2) node {6};
\draw[] (8,-0.2) node {7};

\draw[] (4.5,-0.7) node {$T_{245}$};
 
\end{tikzpicture}     
        
    \end{subfigure}
    \hfill
    \caption{Reduction performed  on the edges $(2,4)$, $(4,5)$ of $T$  is not noncrossing, however when performed on $T_{245}$ it is noncrossing.} 
    \label{fig:tklm} 
\end{figure}
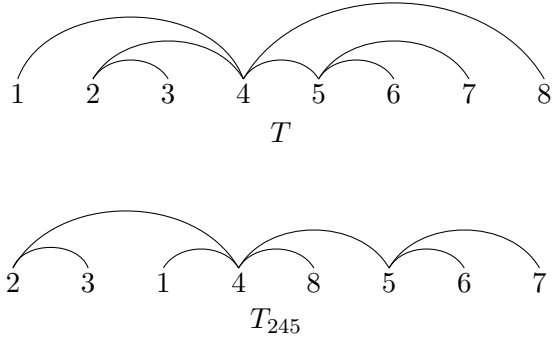

 \bl \label{lem:!!} For a noncrossing tree $T$  on the vertex set $[n]$ and any two edges $(k,l), (l,m)$ of $T$ that are nonalternating, we have that 
 \be \label{eq:noncross} R^{noncross}_T(t_i \mid 1\leq i\leq n-1)=R^{noncross}_{T_{klm}}(t_{i} \mid 1\leq i\leq n-1). \ee
 \el
 
 \proof We prove this lemma by induction on the number of increasing paths in $T$. Suppose there is a vertex $v\neq l$ that is nonalternating. Perform noncrossing reductions at $v$ in both $T$ and $T_{klm}$ obtaining three descendants. Note that the graphs obtained in this fashion from $T$ and $T_{klm}$ are in natural bijection, and they each have fewer number of increasing paths than does $T$, so by inductive hypothesis the Lemma is true for them. However, $R^{noncross}_T(t_i \mid 1\leq i\leq n-1)$ and $R^{noncross}_{T_{klm}}(t_i \mid 1\leq i\leq n-1)$ is the sum of the t-reduced forms corresponding to the mentioned graphs, so we are done. 
 
 It remains to prove the case when the only nonalternating vertex of $T$ is $l$. This is accomplished in Lemma \ref{lem:ll} below. \qed


 \bl \label{lem:ll} For $T:=T^l=([n], \{(i,l), (l, j) \mid 1\leq i<l, l<j\leq n\})$, for some $2\leq l \leq n-1$, and any two edges $(k,l), (l,m)$ of $T$ that are nonalternating, we have that 
 \be \label{eq:noncross} R^{noncross}_T(t_i \mid 1\leq i\leq n-1)=R^{noncross}_{T_{klm}}(t_{i} \mid 1\leq i\leq n-1). \ee
 \el
 
 \proof The only both left and right pseudo-component of $T$ is $\{(1,l), (l, n)\}$, its left pseudo-components are $\{(i,l) \mid 2\leq i\leq l-1\}$, and its right pseudo-components are $\{(l,i) \mid l+1\leq i\leq n-1\}$. Similarly, the only both left and right pseudo-component of $T$ is $\{(k,l), (l, m)\}$, its left pseudo-components are $\{(i,l)\mid i \neq k$, $1\leq i\leq l-1\}$, and its right pseudo-components are $\{(l,i)\mid i\neq m, l+1\leq i\leq n\}$.
 Using this one can prove by induction on $l$ that there is a bijection between the forests described in Theorem \ref{gen:non} for $T$ and for those of $T_{klm}$ such that the number of edges emanating from any vertex $i  \in [n]$ is preserved.  While such a proof is not hard, it is technical to describe, and we leave it to the interested reader.
 \qed

\bigskip

\noindent{\it Proof of Theorem \ref{thm:unique}.}  We proceed by induction on the number of increasing paths in $T$. If we start by a noncrossing reduction $(k,l), (l,m)$ in $T$, then no matter how we reduce the three descendants of $T$ which each have fewer number of increasing paths, we obtain that the t-reduced form  of $m[T]$ we get  is ${\rm Crf}_T({\bf t})$. 

Suppose we start with a reduction $(k,l), (l,m)$ in $T$ that is a crossing reduction. Redraw the tree $T$ as $T_{klm}$. Since we can apply the inductive hypothesis to all three descendants of $T_{klm}$, we get that the  t-reduced form of $m[T]$ obtained this way is ${\rm Crf}_{T_{klm}}({\bf t})$. 

However, by Lemma \ref{lem:!!} ${\rm Crf}_T({\bf t})={\rm Crf}_{T_{klm}}({\bf t})$, thereby proving the theorem.
\qed

\section*{Acknowledgements} The authors are  grateful to Allen Knutson for  many helpful discussions and references about topics related to this research.  The authors also  thank  Ed Swartz  for several  helpful discussions.  Special thanks  go to Sergey Fomin for an extensive conversation about Grothendieck polynomials. We thank Avery St. Dizier and the anonymous referee for their careful reading and suggestions which resulted in the significant improvement of the exposition.
 
\bibliography{biblio-kir}
\bibliographystyle{alpha}

\end{document}